\documentclass[a4paper, 12pt]{article}

\usepackage{cite}
\usepackage{amsfonts, amsmath, amsthm}
\usepackage{mathrsfs}
\usepackage{amssymb}
\usepackage{latexsym}
\usepackage[dvips]{graphicx}
\usepackage{color}

\theoremstyle{plain}
\newtheorem{thm}{Theorem}[section]
\newtheorem{prop}{Proposition}[section]
\newtheorem{lem}{Lemma}[section]

\theoremstyle{definition}

\theoremstyle{definition}
\newtheorem{rmk}{Remark}[section]

\numberwithin{equation}{section}

\newcommand{\Z}{\mathbb{Z}}

\newcommand{\R}{\mathbb{R}}
\newcommand{\C}{\mathbb{C}}

\newcommand{\pa}{\partial}
\newcommand{\eps}{\varepsilon}
\newcommand{\sig}{\sigma}

\DeclareMathOperator{\realpart}{\rm Re}
\DeclareMathOperator{\imagpart}{\rm Im}
\DeclareMathOperator{\supp}{\rm supp}

\begin{document}
\title{
On the decay estimate for small solutions to 
nonlinear Klein-Gordon equations with dissipative structure 
}
\author{
          Yoshinori Nishii\thanks{
              Department of Mathematics, Graduate School of Science,
              Tokyo University of Science.
              1-3 Shinjuku-ku, 
              Tokyo 162-8601, Japan.
              (E-mail: {\tt yoshinori.nishii@rs.tus.ac.jp})
             }
}
\date{\today }
\maketitle


\noindent{\bf Abstract:}\ 
We consider the Cauchy problem for cubic nonlinear Klein-Gordon equations 
in one space dimension. 
We give the $L^p$-decay estimate for the small data solution 
and show that it decays faster than the free solution 
if the cubic nonlinearity has the suitable dissipative structure. 
\\

\\



\section{Introduction and main result}  \label{sec_intro}
We consider the large time behavior of the solution $u=u(t,x)$ to the Cauchy problem 
\begin{align}\label{eq}
\left\{
\begin{aligned}
 &(\Box +1)u=F(u,\pa_t u,\pa_x u), &t>0, \ x\in\R, \\
 &u(0,x)=\eps f(x), \quad \pa_t u(0,x)=\eps g(x),  &x\in\R,
\end{aligned}
\right.
\end{align}
where $\Box =\pa_t^2-\pa_x^2$ with $\pa_t =\pa/\pa t$, $\pa_x=\pa/\pa x$.
We suppose that $f$, $g:\R\to\R$ are compactly-supported $C^{\infty}$-functions, 
$\eps>0$ is a small parameter and 
$F$ is a cubic homogeneous polynomial in $(u,\pa_t u,\pa_x u)$.

To clarify our motivation, let us recall some of the previous results briefly. 
It is well known that cubic nonlinearities in one space dimension causes 
critical situation in the sense of asymptotic behavior of solutions.
According to the result of Yordanov~\cite{Y}, there is a blowing-up solution 
with arbitrary small data for $F=u_t^2 u_x$ 
(see also Proposition~7.8.8 in \cite{Hor}). 
On the other hands, small data global existence (SDGE) holds for 
some classes of the cubic nonlinearities. 
Georgiev-Yordanov~\cite{GY} consider the case $F=u^3$. 
They show SDGE holds and a global solution 
does not behave like a free solution (see also \cite{HN10}). 
In contrast, Katayama~\cite{Kat}, Moriyama~\cite{Mor} give a class of the 
cubic nonlinearities for which SDGE and asymptotically free hold.
According to Delort~\cite{D}, these SDGE results can be summarized 
by the function $K_F:\R\to\C$ which defined by 
\begin{equation*}
K_F(z)=\frac{i}{2\pi} \int_0^{2\pi} 
F(\cos\theta,-\cosh z\sin\theta , \sinh z \sin\theta)e^{-i\theta} d\theta
\end{equation*}
for $z\in\R$. 
It has been shown that SDGE holds for \eqref{eq} and given asymptotic 
profile of global solutions under
\begin{equation}\label{A0}
\realpart[K_F(z)] = 0, \quad z\in\R. 
\tag{${\bf A}_0$}
\end{equation}
This result is extended by Sunagawa~\cite{S3} as following: 
SDGE holds under the condition 
\begin{equation}\label{A}\tag{{\bf A}}
\inf_{z\in\R}\realpart[K_F(z)] \ge 0.
\end{equation}
Moreover, the global solution $u(t,x)$ to \eqref{eq} under \eqref{A} behaves like
\begin{align*}
u(t,x)=&\frac{\frac{1}{\sqrt{t}} \realpart \left[a(x/t) \exp{\left(
	i(t^2 - |x|^2)_{+}^{1/2} + i\Psi_F(x/t)|a(x/t)|^2\mathscr{L}(t,|a(x/t)|^2\Phi_F(x/t))
 \right)}\right]}
{\sqrt{1+2\Phi_F(x/t)|a(x/t)|^2\log t}} \\
&+ O\left(t^{-1/2}(\log t)^{-3/2}\right)
\end{align*}
as $t\to+\infty$ uniformly in $x\in\R$.
Here $i = \sqrt{-1}$, $(\cdot)_{+} = \max\{ \cdot, 0 \}$, 
$a:\R\to\C$ is a suitable function, $\Phi_F(y)$, $\Psi_F(y)$ are given by 
\begin{align*}
\Phi_F(y)=(1-y^2)^{1/2}\realpart{\left[K_F(\tanh^{-1}y)\right]}, \\
\Psi_F(y)=(1-y^2)^{1/2}\imagpart{\left[K_F(\tanh^{-1}y)\right]} 
\end{align*}
for $-1<y<1$ and
$\mathscr{L}(\tau,\varphi)$ is defined by 
\begin{align} \label{def_L}
\mathscr{L}(s, \varphi)
&=
-\int_{1}^{s} \frac{1}{\sig(1+2\varphi\log\sig)} \, d\sig \notag\\
&=
\left\{
\begin{aligned}
 &-\log s, 
  &&{\rm if}\ \varphi = 0, \\
 &-\frac{1}{2\varphi}\log(1+2\varphi\log s), 
  &&{\rm if}\ \varphi \neq 0, \ 1+2\varphi\log s >0.
\end{aligned}
\right.
\end{align}
The typical example which satisfies \eqref{A} but violating \eqref{A0} 
is nonlinear damping $F=-u_t^3$ (see Nakao~\cite{Nakao}, Mochizuki-Motai~\cite{MM}).
If $\realpart [K_F(z_0)] > 0$ holds for some $z_0\in\R$, 
a solution of \eqref{eq} decays faster than a free solution 
in the sense of $L^{\infty}$ in some region. 
Therefore it is fair to say that \eqref{A} (without \eqref{A0}) is dissipative condition.
Indeed, Kim-Sunagawa~\cite{KimS} shows that the global solution 
$u(t)$ to \eqref{eq} decays faster than non-trivial free solutions 
in $L^p$ with $2\le p \le \infty$ 
under the following conditions stronger than \eqref{A}:
If there exists a positive constant $C_{0}$ such that 
\begin{equation}\label{B0}
\realpart{[K_{F}(z)]} \ge C_{0}\cosh^3 z, 
\tag{{\bf B$_0$}}
\end{equation}
for all $z\in\R$, then 
\begin{equation}\label{decay_B0}
\sum_{|\alpha|\le1}\| \pa^{\alpha} u(t) \|_{L^p}
=
O\left( t^{-(1/2-1/p)}(\log t)^{-1/2} \right)
\end{equation}
as $t\to+\infty$, where $\pa = (\pa_t, \pa_x)$.
If there exists a positive constant $C_{1}$ such that 
\begin{equation}\label{B1}
\realpart[K_{F}(z)] \ge C_{1}\cosh z, 
\tag{{\bf B$_1$}}
\end{equation}
for all $z\in\R$, then 
\begin{equation}\label{decay_B1}
\| u(t) \|_{L^p}
=
O\left( t^{-(1/2-1/p)}(\log t)^{-1/2} \right)
\end{equation}
as $t\to+\infty$.
Note that \eqref{B0} or \eqref{B1} implies $\inf_{z\in\R} \realpart{[K_{F}(z)]} >0$.
Thus, it is natural to ask what happens if 
\eqref{A} and $\inf_{z\in\R} \realpart{[K_{F}(z)]} =0$
are satisfied but \eqref{A0} is violated.
To the author's knowledge, there are no previous works 
which address this question except \cite{Mu} 
(see Remark~\ref{rmk_optimal} below). 

Our aim of this paper is to investigate the $L^p$-decay estimate 
under \eqref{A} without \eqref{A0} and \eqref{B0}. 
As a consequence, we reveal the dissipative property of the missing piece of \eqref{A}. 
In order to state our main results, let us introduce the following three conditions:
\begin{itemize}
\item[$({\bf B}_2)$]
There exists a positive constant $C_2$ such that
\begin{equation*}
\realpart{[K_{F}(z)]} \ge \frac{C_2}{\cosh z}
\end{equation*}
for all $z\in\R$.

\item[$({\bf B}_3)$]
There exists a positive constant $C_3$ such that
\begin{equation*}
\realpart{[K_{F}(z)]} \ge \frac{C_3}{\cosh^3 z}
\end{equation*}
for all $z\in\R$.

\item[$({\bf C})$]
There exists $z_0\in\R$ such that
$\realpart{[K_{F}(z_0)]} = 0$ and $\realpart{[K_{F}(z)]} >0$ for $z\in\R\setminus\{ z_0 \}$.
\end{itemize}

We are in a position to state the main results of this paper.
\begin{thm} \label{thm_B1}
Assume the condition \eqref{B1} is satisfied. 
Then there exists a positive constant $C$, not depending on $\eps$, such that 
the global solution $u(t)$ to \eqref{eq} satisfies 
\begin{equation*}
\| \pa u(t) \|_{L^p} 
\le
\left\{
\begin{aligned}
 &C(\log (2+t))^{-1/2} (\log(1+ \log (2+t)))^{1/2},  &&p=2, \\
 &C(1+t)^{-(1/2-1/p)} (\log (2+t))^{-1/p},  &&2<p\le\infty
\end{aligned}
\right.
\end{equation*}
for $t\ge 0$ if $\eps$ is sufficiently small.
\end{thm}
\begin{thm} \label{thm_B2}
Assume the condition $({\bf B}_2)$ is satisfied. 
Then there exists a positive constant $C$, not depending on $\eps$, such that 
the global solution $u(t)$ to \eqref{eq} satisfies 
\begin{equation*}
\| u(t) \|_{L^p} 
\le
\left\{
\begin{aligned}
 &C(\log (2+t))^{-1/2} (\log(1+ \log (2+t)))^{1/2},  &&p=2, \\
 &C(1+t)^{-(1/2-1/p)} (\log (2+t))^{-1/p},  &&2<p\le\infty
\end{aligned}
\right.
\end{equation*}
and
\begin{equation*}
\| \pa u(t) \|_{L^p} 
\le
 C(1+t)^{-(1/2-1/p)} (\log (2+t))^{-1/(2p)}, \ 2\le p\le \infty
\end{equation*}
for $t\ge 0$ if $\eps$ is sufficiently small.
\end{thm}
\begin{thm} \label{thm_B3}
Assume the condition $({\bf B}_3)$ is satisfied. 
Then there exists a positive constant $C$, not depending on $\eps$, such that 
the global solution $u(t)$ to \eqref{eq} satisfies 
\begin{equation*}
\| u(t) \|_{L^p} 
\le
C(1+t)^{-(1/2-1/p)} (\log (2+t))^{-1/(2p)}
\end{equation*}
and
\begin{equation*}
\| \pa u(t) \|_{L^p} 
\le
 C(1+t)^{-(1/2-1/p)} (\log (2+t))^{-1/(3p)}
\end{equation*}
for $2\le p \le \infty$ and $t\ge 0$ if $\eps$ is sufficiently small.
\end{thm}
\begin{thm} \label{thm_C}
Assume the condition $({\bf C})$ is satisfied. 
Then there exists a positive constant $C$, not depending on $\eps$, such that 
the global solution $u(t)$ to \eqref{eq} satisfies 
\begin{equation*}
\sum_{|l|\le1} \| \pa^l u(t) \|_{L^p} 
\le
 C(1+t)^{-(1/2-1/p)} (\log (2+t))^{-1/(2p)}
\end{equation*}
for $2\le p\le \infty$ and $t\ge 0$ if $\eps$ is sufficiently small.
\end{thm}

\begin{rmk}\label{rmk_nonli}
If we assume the condition $({\bf A})$, 
then only one of $({\bf A}_0)$, $({\bf B}_3)$ or $({\bf C})$ is satisfied.
Consequently, $L^p$-decay property of global solutions to \eqref{eq} 
under $({\bf A})$ can be classified as above theorems. 
We give the detail of this fact and the concrete examples of our conditions 
in Section~\ref{sec_nonli}, below.
\end{rmk}

\begin{rmk}\label{rmk_analogue}
Similar results for the semilinear wave equation 
$\Box u = F(\pa u)$ 
in two space dimensions has been obtained in \cite{KMS}, \cite{KMatsS} and \cite{NST} 
for strongly and weakly dissipative structure, respectively. 
Recently, optimality of these results are discussed in \cite{NSS}.  
We also refer the analogues results for cubic nonlinear Scr\"odinger equation 
$i\pa_t u +\frac{1}{2}\pa_x^2 u = N(u, \pa_x u)$ in one space dimension. 
$L^2$-decay for power type nonlinearity 
$N=\lambda |u|^2u$ with $\imagpart{\lambda}<0$, 
which is a typical example of dissipative nonlinearity, 
is studied in 
\cite{HLN, Hoshino, JJL, KitaShim, OgSat, Sat1, Sat3, Sat4, Shim}. 
For more general strongly dissipative nonlinearities, 
$L^2$-decay is obtained in \cite{LiS} 
and weakly dissipative structure is considered in \cite{LNSS3, LNSS4}. 
\end{rmk}

\begin{rmk}\label{rmk_optimal}
According to \cite{MS}, $L^{\infty}$-decay rate of solutions 
to cubic nonlinear Klein-Gordon equation in one space dimension 
is no faster than $O\left( t^{-1/2} (\log t)^{-1/2} \right)$. 
In this sense, for $p=\infty$, 
we can not improve \eqref{decay_B0} under \eqref{B0} 
and \eqref{decay_B1} under \eqref{B1}. 
We next consider the optimality of Theorem~\ref{thm_B1}. 
The special case $F=-u^2\pa_t u$, which is a concrete example of ({\bf B$_1$}), 
is considered in Murakami \cite{Mu}. 
For this case, it holds that 
$\| \pa_t u(t) \|_{L^{\infty}}$ decays like $O\left( t^{-1/2} (\log t)^{-1/2+\eta} \right)$ 
as $t\to\infty$ with $0<\eta<1/4$. 
By using this estimate, we get
\[
\| \pa_t u(t) \|_{L^{p}} = O\left( t^{-(1/2-1/p)} (\log t)^{-1/2+\delta} \right)
\]
as $t\to\infty$ with $0<\delta\ll 1$ for $F=-u^2\pa_t u$. 
This estimate is sharper than Theorem~\ref{thm_B1} for $\pa_t u$. 
We emphasize that we can treat more general nonlinearities 
and not only $\pa_t u$ but also $\pa_x u$ in Theorem~\ref{thm_B1}. 
However, the estimate of $\| \pa u(t) \|_{L^p}$ may not be optimal.
\end{rmk}

The rest of this paper is organized as follows. 
In Section~\ref{sec_nonli}, we introduce where do our conditions come from 
and give concrete examples of nonlinearities of each conditions.
Before we prove the main theorems, we summaries the reduction argument and 
pointwise asymptotics of solution in Section~\ref{sec_prelim}. 
Section~\ref{sec_proof} is devoted to the proofs of 
Theorems~\ref{thm_B1}, \ref{thm_B2}, \ref{thm_B3} and \ref{thm_C}. 

\section{Classification of nonlinearities and examples}  \label{sec_nonli}
In this section, we introduce how to classify the 
nonlinearities satisfying the condition $({\bf A})$. 
As a consequence, we see all of them are treated 
by ({\bf A$_0$}), ({\bf B$_j$}) for $j=0,1,2,3$ and $({\bf C})$.
We also give some example of nonlinearities which satisfies our conditions. 

\subsection{Classification of nonlinearities}  \label{subsec_class}
We first give the explicit representation of $K_{F}(z)$. 
Without loss of generality, $F$ can be written in the form
\begin{align}\label{nonli_gen}
F(u, u_t, u_x)=
&(\gamma_1 u^2 + \gamma_2 u_t^2 + \gamma_3 u_x^2 + \gamma_4 u_t u_x )u \notag\\
&+(\gamma_5 u^2 + \gamma_6 u_t^2 + \gamma_7 u_x^2 )u_t 
+ (\gamma_8 u^2 + \gamma_9 u_t^2 + \gamma_{10} u_x^2 )u_x
\end{align}
with real constants $\gamma_1, \dots, \gamma_{10}$. 
For this $F$, it holds that
\begin{align*}
K_{F}(z)
=& \frac{i}{8}
 \left( 
  3\gamma_1 + \gamma_2\cosh^2 z + \gamma_3\sinh^2 z - \gamma_4\cosh z \sinh z 
 \right) \\
&-\frac{\cosh z}{8} 
 \left( 
  \gamma_5 + 3\gamma_6 \cosh^2 z + 3\gamma_7 \sinh^2 z 
 \right) \\
&+\frac{\sinh z}{8} \left( 
 \gamma_8 + 3\gamma_9 \cosh^2 z + 3\gamma_{10} \sinh^2 z 
 \right).
\end{align*}

To classify the nonlinearities satisfying ({\bf A}), we focus on some cubic polynomials.
We introduce the following elementary lemma, whose proof is skipped.
\begin{lem}\label{lem_poly}
Let $p\in\R[y]$ with $\deg p \le3$.
If $\inf_{y\in(-1,1)}p(y)\ge0$, 
then we have either of the following six assertions:
\begin{itemize}
\item[$(a_0)$] 
$p(y)=0$ for all $y\in[-1,1]$.
%

\item[$(b_0)$] 
$p(y)>0$ for all $y\in[-1,1]$.
%

\item[$(b_1)$] 
$p(y)>0$ for all $y\in(-1,1)$ and 
there exists $c_1>0$ such that either or both of the followings hold:
\begin{itemize}
\item 
$p(y) = (1-y) \left(c_{1} + o(1)\right)$ as $y\to1$ 
and $p'(-1)\neq0$.
\item 
$p(y) = (1+y) \left(c_{1} + o(1)\right)$ as $y\to-1$
and $p'(1)\neq0$.
\end{itemize}
%

\item[$(b_2)$] 
$p(y)>0$ for all $y\in(-1,1)$ and 
there exists $c_2>0$ such that either of the following holds:
\begin{itemize}
\item 
$p(y) = (1-y)^2 \left(c_{2} + o(1)\right)$ as $y\to1$.
\item 
$p(y) = (1+y)^2 \left(c_{2} + o(1)\right)$ as $y\to-1$.
\end{itemize}
%

\item[$(b_3)$] 
Either $p(y) = c_{3} (1-y)^3$ or $p(y) = c_{3} (1+y)^3$ 
with some $c_3>0$.
%

\item[$(c)$] There exist $y_0\in(-1,1)$ and $c_0>0$ such that 
$p(y)>0$ for $y\in(-1,1)\setminus\{ y_0 \}$, 
$p(y) = (y-y_0)^2 \left(c_{0}+ o(1)\right)$ as $y\to y_0$ 
and $p'(-1), p'(1)\neq0$.

\end{itemize}
\end{lem}
We put 
\begin{equation}\label{P_F}
P_{F}(y) = 8(1-y^2)^{3/2} \realpart[K_{F}(\tanh^{-1} y)], 
\end{equation}
then $P_{F}$ is a (at most) cubic polynomial with real coefficients in $y\in(-1,1)$. 
More precisely, we have
\[
P_{F}(y) 
= 
-(\gamma_5 + 3\gamma_6) + (\gamma_8 + 3\gamma_9)y 
+ (\gamma_5 - 3\gamma_7)y^2 + (3\gamma_{10} - \gamma_8)y^3.
\]
Since ({\bf A}) is equivalent to 
\[
\inf_{y\in(-1,1)} P_F (y)\ge0, 
\]
we can apply Lemma~\ref{lem_poly} as $p=P_F$. 
It is clear that ($a_0$) 
is equivalent to ({\bf A$_0$}). 
From the bijectivity of $z\mapsto\tanh z$, we see that ($c$) 
and ({\bf C}) are equivalent. 
If ($b_j$) 
holds for $P_F$, we have 
\[
P_F(y) \ge 8C_j(1-y^2)^j \ge 8C_j(1-y^2)^3, \quad -1\le y \le 1, 
\]
which implies ({\bf B$_3$}) for $j=0,1,2,3$. 
By the above, we deduce that one of ({\bf A$_0$}), ({\bf B$_3$}) or ({\bf C}) holds  
if the nonlinear term satisfies ({\bf A}).
\qed

\subsection{Examples}  \label{subsec_nonli}

We close this section with the concrete examples of nonlinearities 
satisfying our conditions ({\bf B$_j$}) and ({\bf C}). 
It is obvious that ({\bf B$_{j-1}$}) is stronger condition than 
({\bf B$_j$}) for $j=1,2,3$. Thus, to check the assumption is not empty, 
we give the cubic nonlinear term which satisfies ({\bf B$_{j-1}$}) 
but fails to satisfy ({\bf B$_j$}) for $j=1,2,3$.

\begin{itemize}
\item[(1)]
As we mentioned above, the typical examples of the nonlinearity satisfying 
({\bf B$_1$}) is $F(u,\pa_t u, \pa_x u) = -u^2 \pa_t u$. 
Corresponding $K_F(z)$ and $P_F(y)$ are 
$\frac{1}{8}\cosh z$ and $(1+y)(1-y)$, respectively.

\item[(2)]
The concrete example of ({\bf B$_2$}) is $F(u,\pa_t u, \pa_x u) = -u^2(\pa_t u + \pa_x u)$. 
We have 
$P_F(y)=(1-y)^2(1+y)$, which satisfies ($b_2$). 
We can also check ({\bf B$_2$}) directly as
\[
K_F(z) = \frac{1}{8}(\cosh z - \sinh z) = \frac{1}{8}e^{-z} \ge \frac{1}{8\cosh z}.
\]

\item[(3)]
From Lemma~\ref{lem_poly}, the cubic polynomial satisfying ({\bf B$_3$}) is only 
$P_F (y) = C(1\pm y)^3$ with some $C>0$. 
In this case, we have $\realpart[K_F(z)] = C(\cosh z \pm \sinh z)^3$.
Therefore the nonlinearities which ({\bf B$_3$}) is satisfied but ({\bf B$_2$}) is violated 
is only $F(u,\pa_t u, \pa_x u) = -C(\pa_t u \pm \pa_x u)^3$ with $C>0$  
up to the first term of \eqref{nonli_gen}. 
In this case, we have
\[
\realpart[K_F(z)] 
= 
\frac{C}{8}(\cosh z \pm \sinh z)^3 = \frac{C}{8}e^{\pm 3z} \ge \frac{C}{8\cosh^3 z}.
\]

\item[(4)]
The last example is $F(u,\pa_t u, \pa_x u) = -(\pa_x u)^2 \pa_t u$. 
For this $F$, we have $K_F(z) = \cosh z \sinh^2 z$ and it vanishes only at $z=0$. 
Thus ({\bf C}) with $z_0=0$ is satisfied in this case.
\end{itemize}


\section{Preliminary}  \label{sec_prelim}

In this section, we summarize the reduction argument and pointwise estimates 
established in previous works 
(\cite{D}, \cite{DFX}, \cite{S1}, \cite{S2}, \cite{S3}, \cite{MSU} etc.).
As a consequence, we obtain the leading part of the solution $u(t,x)$ to \eqref{eq}, 
which play an important roll in our proof of the main theorems in the next section. 
In what follows, we denote several positive constants by the same letter C, which may
vary from one line to another.
\subsection{Reduction of the equation}
In this subsection, we make some reduction of the problem 
based on \cite{S3} (see also \cite{D}, \cite{DFX}, \cite{S1}, \cite{S2} and \cite{MSU}). 
The argument is divided into three steps. 

\noindent{\bf \underline{Step 1:}}\ 
First we start with the local existence theorem. 
Since $f$ and $g$ are compactly-supported, we can take $B>1$ such that
\begin{align*}
\supp f \cup \supp g \subset \{ x\in\R; |x|\le B \}.
\end{align*}
Let us remember the following local existence result by Delort \cite[Proposition~1.4]{D}.
\begin{prop}
Let $m\ge3$ be an integer. Suppose that $(f,g)\in H^{m}(\R)\times H^{m-1}(\R)$, 
supported in $\{ x\in\R \,;\, |x|\le B \}$ and fix $\tau_0>1+2B$. 
Then, there exists $\eps_0>0$ with the following properties: 
If $\eps\in(0,\eps_0]$, there exist $T \ge (\tau_0^2 - 3B^2)/(2B)$ and unique solution 
$u$ to \eqref{eq} satisfying $u\in C([0,T];H^m(\R))\cap C^1([0,T];H^{m-1}(\R))$. 
Moreover, $u$ is defined on the curve 
$
\left\{ (t,x)\in\R\times\R; 
t\ge0, |x|\le(\tau_0^2 - B^2)/(2B), 
(t+2B)^2-x^2=\tau_0^2\right\}
$ 
and we have 
\[
\sum_{0\le k+l \le m} \int_{|y|\le\frac{\tau_0^2 - B^2}{2B}, } 
\left| 
(\pa_{t}^{k} \pa_{x}^{l} u)\left(\sqrt{\tau_0^2 + y^2} -2B, y\right) 
\right|^2 \,dy
\le C\eps^2.
\]
Furthermore, the solution $u$ has finite propagation speed: 
\begin{equation}\label{fps}
\supp u(t,\cdot) \subset \{ x\in\R; |x|\le t +B \}
\end{equation}
for $0\le t \le T$.
\end{prop}
This proposition allows us to treat the problem 
on the domain 
\begin{equation*}
\left\{ (t,x)\in\R\times\R  ; (t+2B)^2-|x|^2\ge\tau_0^2, \ t\ge0 \right\}.
\end{equation*}

\noindent{\bf \underline{Step 2:}}\ 
Next, we introduce the hyperbolic coordinates 
\[
\tau = \sqrt{(t+2B)^2 - x^2}, \quad z=\tanh^{-1}\left( \frac{x}{t+2B} \right)
\]
for $|x|< t+2B$.
Note that it is equivalent to
\begin{align}\label{hyp_coord}
t+2B=\tau \cosh z, \quad x=\tau\sinh z, \quad \tau\ge\tau_0, \quad z\in\R
\end{align}
and holds that 
\begin{align} \label{hyp_derivative}
\pa_t = (\cosh z) \pa_{\tau} - \frac{\sinh z}{\tau}\pa_z, \quad
\pa_x = -(\sinh z)\pa_{\tau} + \frac{\cosh z}{\tau}\pa_z, \quad
\Box = \pa_{\tau}^2+ \frac{1}{\tau}\pa_{\tau} - \frac{1}{\tau^2}\pa_{z}^2.
\end{align}
For the solution $u$ to \eqref{eq}, we define the new unknown function 
$v=v(\tau,z)$ by
\begin{align}\label{u_to_v}
u(t,x)=\frac{\chi(z)}{\tau^{1/2}}v(\tau,z), 
\end{align}
where $\chi\in C^{\infty}(\R)$ is a weight function satisfying 
\begin{equation*}
0<\chi(z)\le c_0 e^{-\kappa|z|}, 
\quad 
|\chi^{(j)}(z)|\le c_j \chi(z)
\end{equation*}
with a large parameter $\kappa\gg1$ 
and positive constants $c_j$, $j\in\Z_{\ge0}$. 
Then, it follows from \eqref{hyp_derivative} that $v$ satisfies 
\begin{equation}\label{eq_v}
\left\{
\begin{aligned}
 &Pv = G(\tau, z, v, \pa_{\tau}v, \pa_{z}v), &\tau>\tau_0, \ z\in\R, \\
 &(v, \pa_{\tau}v)|_{\tau=\tau_0} = (\eps v_0, \eps v_1),  &z\in\R,
\end{aligned}
\right.
\end{equation}
where $v_0$, $v_1$ are with compactly-supported smooth functions of $z$, 
\begin{align*}
P =\pa_{\tau}^2 - \frac{1}{\tau^2}
\left( \pa_z^2 + 2\frac{\chi'(z)}{\chi(z)}\pa_z + \frac{\chi''(z)}{\chi(z)} -\frac{1}{4} \right)
 + 1
\end{align*}
and
\begin{align*}
G=\frac{\chi(z)^2}{\tau}F(v, (\cosh z) \pa_{\tau}v, -(\sinh z) \pa_z v) + Q.
\end{align*}
Here $Q$ takes the form of
\begin{equation*}
Q=
\sum_{\nu_1, \nu_2 =0}^1 \,
\sum_{\substack{0\le\mu_1\le\nu_1 \\ 0\le\mu_2\le\nu_2 \\ 0\le\mu_3\le1}}
\frac{q_{\nu_1, \nu_2}^{\mu_1,\mu_2,\mu_3}(z)}{\tau^{2+\nu_1+\nu_2}}
\left( \pa_{\tau}^{1-\nu_1}\pa_z^{\mu_1}v \right)
\left( \pa_{\tau}^{1-\nu_2}\pa_z^{\mu_2}v \right)
\left( \pa_z^{\mu_3}v \right)
\end{equation*}
with some $q_{\nu_1, \nu_2}^{\mu_1,\mu_2,\mu_3}\in C^{\infty}(\R)$ satisfying 
\begin{equation*}
|\pa_{z}^j q_{\nu_1, \nu_2}^{\mu_1,\mu_2,\mu_3}(z)| \le c_j e^{-(2\kappa-3)|z|}
\end{equation*}
for $j\in\Z_{\ge0}$.
We define 
\[
 E_s(\tau)
=
 \sum_{k=0}^{s} \frac{1}{2} \int_{\R} 
 \left(
  \left| \pa_{\tau}\pa_{z}^k v(\tau,z) \right|^2
  + \left| \frac{\pa_{z}}{\tau}\pa_{z}^k v(\tau,z) \right|^2 
  + \left| \pa_{z}^k v(\tau,z) \right|^2 
 \right)
  \, dz
\]
for $s\in\Z_{\ge0}$, 
then we have the following estimates.
\begin{lem}\label{lem_apriori_v}
Let $\eps$ be suitably small. Assume that \eqref{A} is satisfied. 
Then, the solution $v$ to \eqref{eq_v} exists globally in $[\tau_0,\infty)$ and 
satisfies
\begin{align}\label{apriori_v}
|v(\tau,z)| + |\pa_{\tau}v(\tau,z)| + \frac{1}{\tau}|\pa_z v(\tau,z)|
\le C\eps
\end{align}
for $(\tau,z)\in[\tau_0,\infty)\times\R$ 
and 
\begin{align}\label{apriori_E}
E_s(\tau)
\le 
C\eps^2\tau^{\delta}
\end{align}
for $s\ge3$ and $0<\delta<1/3$
for $\tau\ge\tau_0$.
\end{lem}
For the proof, see \cite{S3} (see also \cite{KimS} and \cite{S2}).
From \eqref{apriori_E}, we also have 
\begin{align}\label{apriori_v_z}
 |\pa_z^j v(\tau,z)|
\le
 C \| v(\tau) \|_{H^{j+1}(\R)}
\le
 C E_{j+1}(\tau)^{1/2}
\le 
C\eps\tau^{\delta/2}
\end{align}
for $\tau\ge\tau_0$, $z\in\R$, $j=1,2$.

\noindent{\bf \underline{Step 3:}}\ 
We make one more reduction. 
Let us define the function $\alpha:[\tau_0,\infty)\times\R\to\C$ by 
\begin{align}\label{v_to_alpha}
\alpha(\tau,z) = e^{-i\tau}\left( 1+\frac{1}{i}\frac{\pa}{\pa \tau} \right) v(\tau,z).
\end{align}
We have already known the following.
\begin{lem}\label{lem_profile}
Under the assumption of Lemma~\ref{lem_apriori_v}, 
we have $\alpha(\tau,z)=\beta(\tau,z)+S(\tau,z)$ with 
\[
|S(\tau,z)|\le \frac{C\eps^3}{\tau}
\]
for $\tau\ge\tau_0$, $z\in\R$ and $\beta(\tau,z)$ satisfies 
\begin{align}\label{eq_beta}
\frac{\pa \beta}{\pa \tau} 
= -\frac{\chi(z)^2 K_F(z)}{\tau}|\beta|^2\beta 
 + R(\tau,z),
\end{align}
\begin{align*}
|R(\tau,z)| \le C\eps\tau^{-2+3\delta/2}
\end{align*}
for $\tau\ge\tau_0$, $z\in\R$.
\end{lem}
The proof can be seen in \cite{S3} (see also \cite{KimS}). 
As we see in the next subsection, the solution $\beta(\tau,z)$ to \eqref{eq_beta} 
gives the leading part of $u(t,x)$ to \eqref{eq}. 
Therefore we call \eqref{eq_beta} the profile equation.

\subsection{Pointwise asymptotics}
In this subsection, we obtain the leading part of $u$, $\pa_t u$ and $\pa_x u$
through the profile equation \eqref{eq_beta}. 
In order to focus on large-time asymptotics of solutions to \eqref{eq_beta}, 
we introduce the following lemma due to Sunagawa \cite[Lemma~2.1]{S3}. 
\begin{lem}\label{lem_ode}
(1)\ 
Let $\beta_0\in\C$ and $\kappa\in\C$ with $\realpart\kappa\ge0$. 
Let $\rho(\tau)$ be a complex-valued function of $\tau\in \R$ which satisfies 
\[
|\rho(\tau)|\le C\tau^{-1-\lambda}, 
\qquad 
\tau\ge\tau_0
\]
with some constants $\lambda>0$ and $\tau_0>1$. 
If $\beta(\tau)$ solves the ODE 
\begin{align*}
\frac{d\beta}{d\tau}
=
-\frac{\kappa}{\tau}|\beta|^2\beta + \eps\rho(\tau), 
\quad
\beta(\tau_0)=\eps\beta_0
\end{align*}
for sufficiently small $\eps$, 
then there exists $\beta_{\infty}\in\C$ such that 
\begin{align}\label{asympt_ode}
 \beta(\tau) 
=
 \frac{\beta_{\infty} \exp\left(i \imagpart{\kappa} |\beta_{\infty}|^2 
         \mathscr{L}(\tau, \realpart{\kappa} |\beta_{\infty}|^2)\right) }
        {\sqrt{1+2\realpart{\kappa} |\beta_{\infty}|^2 \log\tau}}
 +
 O\left(\left( \log\tau \right)^{-3/2}\right)
\end{align}
as $\tau\to+\infty$, where $\mathscr{L}$ is given by \eqref{def_L}.

\noindent (2)\ 
If $\beta(\tau)$, $\beta_0$, $\kappa$, $\rho(\tau)$ depend smoothly 
on some parameter $z\in\R$ and satisfy 
\begin{equation*}
|\pa_z^j \beta(\tau,z)|\le M_j\tau^{\delta}, 
\quad 
|\pa_z^j \kappa(z)|\le M_j, 
\quad 
|\pa_z^j \rho(\tau,z)| \le M_j \tau^{-1-\lambda_j}
\end{equation*}
for $j\in\Z_{\ge0}$ 
with some constants $M_j>0$, $\lambda_j>0$ and sufficiently small $\delta>0$, 
then $\beta_{\infty}$ is a bounded smooth function of $z$ 
and \eqref{asympt_ode} is valid uniformly with respect to $z$.

\noindent (3)\ 
If, in addition, $\kappa(z)$ decays like $O(1/z)$ as $|z|\to\infty$, 
then \eqref{asympt_ode} can be replaced by
\begin{align*}
 \beta(\tau,z) 
=&
 \frac{\beta_{\infty}(z) \exp\left(i \imagpart{\kappa(z)} |\beta_{\infty}(z)|^2 
         \mathscr{L}(\log(\tau\cosh z), \realpart{\kappa(z)} |\beta_{\infty}(z)|^2)\right) }
        {\sqrt{1+2\realpart{\kappa(z)} |\beta_{\infty}(z)|^2 \log(\tau\cosh z)}} \\
 &+
 O\left(e^{h|z|}\left( \log(\tau\cosh z) \right)^{-3/2}\right)
\end{align*}
with arbitrary small positive number $h$, 

\end{lem}
For the proof, see \cite{S3}. 
We apply this lemma with 
$\kappa(z) = \chi(z)^2 K_F(z)$ and 
$\rho(\tau,z) = R(\tau,z)$ appears in Lemma~\ref{lem_profile},  
then the leading part of $\beta(\tau,z)$ as $\tau\to\infty$ is given by 
\begin{align*}
 A(\tau,z) 
:=
 \frac{\beta_{\infty}(z) \exp\left(i \chi(z)^2 \imagpart{[K_F(z)]} |\beta_{\infty}(z)|^2 
         \mathscr{L}(\log(\tau\cosh z), \chi(z)^2\realpart{[K_F(z)]} |\beta_{\infty}(z)|^2)\right) }
        {\sqrt{1+2\chi(z)^2\realpart{[K_F(z)]} |\beta_{\infty}(z)|^2 \log(\tau\cosh z)}} 
\end{align*}
uniformly in $z\in\R$.
Since $\beta_{\infty}$ is a bounded function of $z$, we have 
\begin{align}\label{A_bdd}
|A(\tau,z)|
=
 \frac{|\beta_{\infty}(z)|}
        {\sqrt{1+2\chi(z)^2\realpart{[K_F(z)]} |\beta_{\infty}(z)|^2 \log(\tau\cosh z)}} 
\le C 
\end{align}
for $\tau\ge\tau_0$, $z\in\R$.
We also have 
\begin{align}\label{beta_to_A}
|\beta(\tau,z) - A(\tau,z)| 
\le 
 \frac{C e^{h|z|}}{\left( \log(\tau\cosh z) \right)^{3/2}}
\end{align}
with $0<h\ll1$ for $\tau\ge\tau_0$, $z\in\R$.
By calculating \eqref{u_to_v} and \eqref{v_to_alpha} with \eqref{hyp_derivative}, 
we obtain 
\begin{align}\label{u}
u(t,x)
=&
\frac{\chi(z)}{\tau^{1/2}} \realpart[\alpha(\tau,z)e^{i\tau}], 
\end{align}
\begin{align}\label{pa_t_u}
\pa_t u(t,x)
=&
 \left( (\cosh z)\pa_{\tau} - \frac{\sinh z}{\tau}\pa_z \right) 
 \left( \frac{\chi(z)}{\tau^{1/2}} v(\tau,z) \right) \notag\\
=
&-\frac{\chi(z) \cosh z }{\tau^{1/2}} \imagpart[\alpha(\tau,z)e^{i\tau}] \notag\\
&-\frac{\chi(z) \cosh z }{\tau^{3/2}} 
  \left( (\tanh z)\pa_{z} + (\tanh z)\frac{\chi'(z)}{\chi(z)} + \frac{1}{2} \right) v(\tau,z) 
\end{align}
and
\begin{align}\label{pa_x_u}
\pa_x u(t,x)
=&
 \left( -(\sinh z)\pa_{\tau} + \frac{\cosh z}{\tau}\pa_z \right) 
 \left( \frac{\chi(z)}{\tau^{1/2}} v(\tau,z) \right) \notag\\
=
&\frac{\chi(z) \tanh z \cosh z }{\tau^{1/2}} \imagpart[\alpha(\tau,z)e^{i\tau}] \notag\\
&+\frac{\chi(z) \cosh z }{\tau^{3/2}} 
  \left( \pa_{z} + \frac{\chi'(z)}{\chi(z)} + \frac{\tanh z}{2} \right) v(\tau,z)
\end{align}
with $\tau=\sqrt{(t+2B)^2 - x^2}$ and $z=\tanh^{-1}\frac{x}{t+2B}$ for $|x|<t+2B$.
By using Lemma~\ref{lem_profile} and \eqref{beta_to_A} for \eqref{u}, 
we get
\begin{align}\label{leading_u}
 |u(t,x)|
&\le
 \frac{\chi(z)}{\tau^{1/2}} 
  \left(
   |A(\tau,z)| + \frac{C e^{h|z|}}{(\log(\tau\cosh z))^{3/2}} + \frac{C\eps^3}{\tau}
  \right) \notag\\
&\le
 \frac{\chi(z)(\cosh z)^{1/2}}{(\tau\cosh z)^{1/2}} |A(\tau,z)| 
 +\frac{C e^{-(\eta -h-1/2)|z|}}{(\tau\cosh z)^{1/2} (\log(\tau\cosh z))^{3/2}}. 
\end{align}
Similarly, by using Lemma~\ref{lem_profile}, \eqref{beta_to_A}, 
\eqref{apriori_v} and \eqref{apriori_v_z} for \eqref{pa_t_u} and \eqref{pa_x_u}, 
we get
\begin{align}\label{leading_pa_u}
 |\pa u(t,x)|
&\le 
 \frac{\chi(z) \cosh z }{\tau^{1/2}} 
  \left(
   |A(\tau,z)| + \frac{C e^{h|z|}}{(\log(\tau\cosh z))^{3/2}} + \frac{C\eps^3}{\tau}
  \right) 
 +\frac{\chi(z) \cosh z }{\tau^{3/2}} 
  \left(
    C\eps \tau^{\delta/2}  + C\eps
  \right) \notag\\
&\le 
 \frac{\chi(z) (\cosh z)^{3/2} }{(\tau\cosh z)^{1/2}} |A(\tau,z)| 
 +\frac{C e^{-(\eta -h-3/2)|z|}}{(\tau\cosh z)^{1/2} (\log(\tau\cosh z))^{3/2}}
\end{align}
with $\tau=\sqrt{(t+2B)^2 - x^2}$ and $z=\tanh^{-1}\frac{x}{t+2B}$ for $|x|<t+2B$.

\section{Proofs}  \label{sec_proof}
In this section, we prove the main theorems. 
We remark that the case $p=\infty$ of 
Theorems~\ref{thm_B1}, \ref{thm_B2}, \ref{thm_B3} and \ref{thm_C} 
are follow from \eqref{leading_u}, \eqref{leading_pa_u} and \eqref{A_bdd}. 
Therefore, we show the case $2\le p<\infty$ in the followings. 
To simplify the calculations, we put 
\[
I^{p}_{m} (t)
:=
 \int_{0}^{2} 
 \left(
\frac{1}
{1+ \eta^{m} \log(t+2)} \\
 \right)^{p/2}
 \,d\eta
\]
for $m\ge1$ and $2\le p <\infty$.
Then we can easily check that  
\begin{align}
I^{p}_{m} (t)
&\le
 C\int_{0}^{2} 
  \frac{1}{1+ \eta^{mp/2} (\log(2+t))^{p/2}} 
 \,d\eta \notag \\
&\le
 C(\log(t+2))^{-1/m}
  \int_{0}^{C(\log(2+t))^{1/m}} 
  \frac{1}{1+ \xi^{mp/2}} 
 \,d\xi \notag \\
&\le 
\left\{
\begin{aligned}\label{I}
 &C(\log(2+t))^{-1} \log(1+\log(t+2)), &&m=1,p=2, \\
 &C(\log(2+t))^{-1/m},  &&mp>2.
\end{aligned}
\right.
\end{align}
\subsection{Proof of Theorems~\ref{thm_B1}, \ref{thm_B2} and \ref{thm_B3}}
  \label{subsec_proof}

Note that under the condition ({\bf B$_j$}), we have 
\[
\frac{\realpart{K_{F}(z)}}{(\cosh z)^{2k+1}}
\ge 
 C_j\left( \frac{1}{\cosh^2 z} \right)^{j+k-1}
=
 C_j (1- \tanh^2 z)^{j+k-1}, 
\]
whence  
\begin{align}\label{u_to_A}
\chi(z)^2 (\cosh z)^{2k+1} |A(\tau,z)|^2
&=
\frac{\left(\chi(z) (\cosh z)^{k+1/2}|\beta_{\infty}(z)|\right)^2 }
{1+2\left(\chi(z) (\cosh z)^{k+1/2}|\beta_{\infty}(z)|\right)^2 
\dfrac{\realpart{K_{F}(z)}}{(\cosh z)^{2k+1}} \log(\tau\cosh z)} \notag\\
&\le
\frac{C}
{1+  C (1- \tanh^2 z)^{j+k-1} \log(\tau\cosh z)} 
\end{align}
for $j=0,1,2,3$ and $k=0,1$.
First, we show the estimate for $u(t)$. 
We assume ({\bf B$_2$}). 
Then, from \eqref{leading_u}, \eqref{u_to_A} with $j=2$, $k=0$ 
and \eqref{hyp_coord}, we have
\begin{align}\label{u_B2}
 |u(t,x)|
&\le 
 \frac{C}{(\tau\cosh z)^{1/2}} \left(
  \left(
   \chi(z)^2 \cosh z |A(\tau,z)|^2
  \right)^{1/2} 
 +\frac{e^{-(\eta -h-1/2)|z|}}{(\log(\tau\cosh z))^{3/2}} \right)\notag\\
&\le 
 \frac{C}{(\tau\cosh z)^{1/2}} 
  \left(
   \frac{1}{1+ C (1- \tanh^2 z) \log(\tau\cosh z)} 
  \right)^{1/2} \notag\\
&\le 
 \frac{C}{(t+2B)^{1/2}} 
  \left(\frac{1}{1+ \left(1- \left( \frac{x}{t+2B}\right)^2 \right) \log(t+2B)} 
  \right)^{1/2}.  
\end{align}
By using the finite propagation speed \eqref{fps}, \eqref{u_B2} and \eqref{I}, 
we obtain 
\begin{align*}
 \| u(t)\|_{L^p}^p
&\le 
 \frac{C}{(t+2B)^{p/2}} 
 \int_{|x|\le t+B} 
   \left(
   \frac{1}{1+ (1-(\frac{x}{t+2B})^2)^{j-1} \log(t+2B)} 
   \right)^{p/2}
 \, dx \notag \\
&\le 
 C(1+t)^{-p/2+1} 
 \int_{-1}^{1} 
   \left(
   \frac{1}{1+ (1-y^2)^{j-1} \log(2+t)} 
   \right)^{p/2}
 \, dy \notag \\
&\le 
 C(1+t)^{-p/2+1} 
 \int_{0}^{1} 
   \left(
   \frac{1}{1+ \eta \log(2+t)} 
   \right)^{p/2}
 \, d\eta \notag \\
&\le 
 C(1+t)^{-p/2+1} I_{1}^{p}(t), 
\end{align*}
which yields
\[
\| u(t) \|_{L^p} 
\le
\left\{
\begin{aligned}
 &C(\log (2+t))^{-1/2} (\log(1+ \log (2+t)))^{1/2},  &&p=2, \\
 &C(1+t)^{-(1/2-1/p)} (\log (2+t))^{-1/p},  &&2<p<\infty
\end{aligned}
\right.
\]
as desired. 
In the same way, we can get 
\[
 \| u(t)\|_{L^p}^p
\le 
 C(1+t)^{-p/2+1} I_{2}^{p}(t) 
\]
and 
\[
\| u(t) \|_{L^p} 
\le
C(1+t)^{-(1/2-1/p)} (\log (2+t))^{-1/(2p)}
\]
for ({\bf B$_3$}) with $2\le p< \infty$. 
Then we finish to prove the first half of Theorems~\ref{thm_B2} and \ref{thm_B3}. 
Similar argument works for the proof of Theorem~\ref{thm_B1} and 
the latter half of Theorems~\ref{thm_B2} and \ref{thm_B3}. 
We assume ({\bf B$_j$}) for $j=1,2,3$.  
Then, by \eqref{leading_pa_u}, \eqref{u_to_A} with $k=1$ 
and \eqref{hyp_coord}, we have
\begin{align*}
 |\pa u(t,x)|
&\le 
 \frac{C}{(\tau\cosh z)^{1/2}} \left(
  \left(
   \chi(z)^2 \cosh^3 z |A(\tau,z)|^2
  \right)^{1/2} 
 +\frac{e^{-(\eta -h-3/2)|z|}}{(\log(\tau\cosh z))^{3/2}} \right)\\
&\le 
 \frac{C}{(\tau\cosh z)^{1/2}} 
  \left(\frac{1}{1+ C (1- \tanh^2 z)^{j} \log(\tau\cosh z)} 
  \right)^{1/2} \\
&\le 
 \frac{C}{(t+2B)^{1/2}} 
  \left(\frac{1}{1+ C \left(1- \left( \frac{x}{t+2B}\right)^2 \right)^{j} \log(t+2B)} 
  \right)^{1/2} 
\end{align*}
and it follows from \eqref{fps} and \eqref{I} that 
we obtain 
\begin{align*}
 \| \pa u(t)\|_{L^p(\R)}^p
&\le 
 \frac{C}{(t+2B)^{p/2}} 
 \int_{|x|\le t+B} 
   \left(
   \frac{1}{1+ (1-(\frac{x}{t+2B})^2)^{j} \log(t+2B)} 
   \right)^{p/2}
 \, dx \notag \\
&\le 
 C(1+t)^{-p/2+1} 
 \int_{-1}^{1} 
   \left(
   \frac{1}{1+ (1-y^2)^{j} \log(2+t)} 
   \right)^{p/2}
 \, dy \notag \\
&\le 
 C(1+t)^{-p/2+1} 
 \int_{0}^{1} 
   \left(
   \frac{1}{1+ \eta^{j} \log(2+t)} 
   \right)^{p/2}
 \, d\eta \notag \\
&\le 
 C(1+t)^{-p/2+1} I_{j}^{p}(t)
\end{align*}
for $j=1,2,3$. 
This completes the proof of Theorem~\ref{thm_B1}, \ref{thm_B2} and \ref{thm_B3}. 
\qed
\subsection{Proof of Theorem~\ref{thm_C}}  \label{subsec_proof_C}
Finally, we prove Theorem~\ref{thm_C}. We assume the condition ({\bf C}). 
Then, by the argument in Subsection~\ref{subsec_nonli}, 
($c$) in Lemma~\ref{lem_poly} with $y_0 = \tanh z_0$ holds for $P_F(y)$, 
namely there exist $c_0>0$ such that 
$P_F(y)>0$ for $y\in(-1,1)\setminus\{ y_0 \}$, 
$P_F(y) = (y-y_0)^2 \left(c_{0}+ o(1)\right)$ as $y\to y_0$ 
and $P_F'(-1), P_F'(1)\neq0$. 
We can take $\delta>0$ so small that there exists $c_{*}>0$ such that 
\begin{align}\label{divide}
P_F(y)\ge
\left\{
\begin{aligned}
 &c_{*}(1+y), &&y\in (-1, y_0 -\delta], \\
 &\frac{c_0}{2}(y-y_0)^2, &&y\in (y_0 -\delta, y_0+\delta), \\
 &c_{*}(1-y), &&y\in [y_0+\delta,1).
\end{aligned}
\right.
\end{align}
By \eqref{P_F}, we also have 
\begin{align}\label{A_to_P_F}
\chi(z)^2 \cosh^3 z |A(\tau,z)|^2
&=
\frac{\left(\chi(z) (\cosh z)^{3/2}|\beta_{\infty}(z)|\right)^2 }
{1+2\left(\chi(z) (\cosh z)^{3/2}|\beta_{\infty}(z)|\right)^2 
\dfrac{\realpart{[K_{F}(z)]}}{(\cosh z)^{3}} \log(\tau\cosh z)} \notag\\
&\le
\frac{C}
{1+  C P_F(\tanh z) \log(\tau\cosh z)}. 
\end{align}
Therefore, 
it follows from \eqref{leading_u}, \eqref{leading_pa_u}, \eqref{A_to_P_F} 
and \eqref{hyp_coord} that 
\begin{align*}
\sum_{|l|\le1} |\pa^{l} u(t,x) |
&\le 
 \frac{C}{(\tau\cosh z)^{1/2}} \left(
  \left(
   \chi(z)^2 \cosh^3 z |A(\tau,z)|^2
  \right)^{1/2} 
 +\frac{e^{-(\eta -h-3/2)|z|}}{(\log(\tau\cosh z))^{3/2}} \right)\\
&\le 
 \frac{C}{(\tau\cosh z)^{1/2}}
  \left(
 \frac{1}{1+ C P_F(\tanh z) \log(\tau\cosh z)} 
  \right)^{1/2} \\
&\le 
 \frac{C}{(t+2B)^{1/2}}
  \left(
 \frac{1}{1+ C P_F(\frac{x}{t+2B}) \log(t+2B)} 
  \right)^{1/2}. 
\end{align*}
By \eqref{fps}, we have 
\begin{align*}
\sum_{|l|\le1} \|\pa^{l} u(t)\|_{L^p}^p
%
%
&\le 
 \frac{C}{(t+2B)^{p/2}} 
 \int_{|x|\le t+B} 
  \left(\frac{1}{1+ C P_F(\frac{x}{t+2B}) \log(t+2B)} 
  \right)^{p/2} 
 \, dx \notag \\
&\le 
 C(1+t)^{-p/2+1} 
 \int_{-1}^{1} 
  \left(\frac{1}{1+ C P_F(y) \log(2+t)} 
  \right)^{p/2} 
 \, dy.
\end{align*}
Now we use \eqref{divide} and \eqref{I}, we get 
\begin{align*}
\sum_{|l|\le1} \|\pa^{l} u(t)\|_{L^p(\R)}^p
&\le 
 C(1+t)^{-p/2+1} 
\Biggl(
 \int_{-1}^{y_0 -\delta} 
  \left(\frac{1}{1+  (1+y) \log(2+t)} 
  \right)^{p/2} 
 dy \notag \\
&\qquad \qquad \qquad \qquad\ +\int_{y_0 -\delta}^{y_0+\delta}
  \left(\frac{1}{1+  (y-y_0)^2 \log(2+t)} 
  \right)^{p/2} 
 dy \notag\\
&\qquad \qquad \qquad\qquad \ +
 \int_{y_0+\delta}^{1} 
  \left(\frac{1}{1+  (1-y) \log(2+t)} 
  \right)^{p/2} 
 dy 
\Biggr) \\
&\le 
 C(1+t)^{-p/2+1} (I^p_{1}(t) + I^p_{2}(t)+I^p_{1}(t) ) \\
&\le 
 C(1+t)^{-p/2+1}(\log(t+2))^{-1/2}.
\end{align*}
Finally, we obtain
\begin{equation*}
\sum_{|l|\le1} \| \pa^l u(t) \|_{L^p} 
\le
 C(1+t)^{-(1/2-1/p)} (\log (2+t))^{-1/(2p)}
\end{equation*}
for $2\le p< \infty$ as desired.
\qed



\medskip
\subsection*{Acknowledgments}
The author would like to thank Professors Hideaki Sunagawa and Satoshi Masaki 
for their useful conversations on this subject.
This work was supported by JSPS Grant-in-Aid for JSPS Fellows (22KJ2801). 
This work was partly supported by MEXT Promotion of Distinctive Joint Research Center Program JPMXP0723833165.


\end{document}